\documentclass[12pt]{article}
\usepackage{amsmath,amsfonts,amsthm,amssymb,eucal}
\begin{document}
\newtheorem{lem}{Lemma}
\newtheorem{teo}{Theorem}
\newtheorem{prop}{Proposition}
\pagestyle{plain}
\title{$p$-Adic Spherical Coordinates and Their Applications}
\author{Anatoly N. Kochubei\footnote{Partially supported by
the Ukrainian Foundation
for Fundamental Research, Grant 29.1/003.}\\
\footnotesize Institute of Mathematics,\\
\footnotesize National Academy of Sciences of Ukraine,\\
\footnotesize Tereshchenkivska 3, Kiev, 01601 Ukraine
\\ \footnotesize E-mail: \ kochubei@i.com.ua}
\date{}
\maketitle

\bigskip
\begin{abstract}
On the space $\mathbb Q_p^n$, where $p\ne 2$ and $p$ does not divide $n$, we construct a $p$-adic counterpart of spherical coordinates. As applications, a description of homogeneous distributions on $\mathbb Q_p^n$ and a skew product decomposition of $p$-adic L\'evy processes are given.
\end{abstract}

\bigskip
{\bf Key words: }\ $p$-adic numbers, spherical coordinates, homogeneous distributions, L\'evy processes

{\bf MSC 2000}. Primary: 11S80. Secondary: 46F05; 60G51.

\bigskip
\section{INTRODUCTION}

Spherical coordinates in $\mathbb R^n$ are among the basic tools of real analysis from its very early days. The usefulness of the decomposition $\mathbb R^n\setminus \{0\}=S^{n-1}\times \mathbb R_+$ is due to the fact that both factors on the right are {\it smooth} manifolds of dimensions smaller than $n$. Thus various $n$-dimensional objects of analysis and geometry are reduced to objects of similar nature in smaller dimensions.

A straightforward generalization to the case of the $n$-dimensional $p$-adic space $\mathbb Q_p^n$ leads to a different situation. What is usually called a $p$-adic unit sphere, the set
$$
\mathcal S_1=\left\{ x=(x_1,\ldots ,x_n)\in \mathbb Q_p^n :\ \max\limits_{1\le j\le n}|x_j|_p=1\right\} ,
$$
is actually an open-closed subset of $\mathbb Q_p^n$, so that it has the same dimension as the ambient space $\mathbb Q_p^n$. The ``radial'' component in the decomposition $\mathbb Q_p^n \setminus \{ 0\}=\mathcal S_1\times p^\mathbb Z$, $p^\mathbb Z=\{ p^N,\ N\in \mathbb Z\}$, given by the equality $x=\left\{ (x_1,\ldots ,x_n)p^N\right\} p^{-N}$, where $\max\limits_{1\le j\le n}|x_j|_p=p^N$, is discrete and does not have the same nature as the ``spherical'' component.

In this paper, assuming that $p\ne 2$ and $p$ does not divide $n$, we construct a coordinate system in $\mathbb Q_p^n$ resembling the classical spherical coordinates. The idea is to identify $\mathbb Q_p^n$ with the unramified extension $K$ of the field $\mathbb Q_p$ of degree $n$ and to use such related objects as the Frobenius automorphism and the norm map. The counterpart of the sphere introduced below is a direct product of a finite set by a hypersurface of the group of principal units of $K$; the counterpart of $\mathbb R_+$ is a multiplicative subgroup of $\mathbb Q_p$ generated by $p^\mathbb Z$ and (an interesting coincidence of terminology!) the group of positive elements of $\mathbb Q_p$ \cite{Sch}. For $n=2$, our construction is different from the polar coordinates introduced in \cite{GGP,VVZ} though the constructions have some common features.

As applications, we obtain, following \cite{Le}, a description of all homogeneous distributions on $\mathbb Q_p^n$ (earlier such a result was known only for $n=1$; only an example was considered for an arbitrary $n$ in \cite{VVZ}), and a skew product representation for $p$-adic L\'evy processes. The latter result follows a recent work by Liao \cite{L} who considered a decomposition of a Markov process on a manifold invariant under a Lie group action; for earlier results regarding decompositions of a Brownian motion into a skew product of the radial motion and the spherical Brownian motion with a time change see \cite{G,PR}.

\section{Preliminaries}

Let us recall some notions and results from $p$-adic analysis and algebraic number theory, which will be used in a sequel. Note that elementary notions and facts regarding $p$-adic numbers and their properties are used without explanations; see \cite{VVZ}. For further details see \cite{FV,Koch,K,R,Sch,W}.

Let $p$ be a prime number, $p\ne 2$. A field $K$ is called a finite extension of degree $n$ of the field $\mathbb Q_p$ of $p$-adic numbers, if $\mathbb Q_p$ is a subfield of $K$, and $K$ is a finite-dimensional vector space over $\mathbb Q_p$, with $\dim K=n$.

For each element $x\in K$, consider a $\mathbb Q_p$-linear operator $L_x$ on $K$ defined as $L_xz=xz$, $z\in K$ (the multiplication in $K$). Its determinant $N(x)=\det L_x$ is an element of $\mathbb Q_p$. The mapping $x\mapsto N(x)$, $K\to \mathbb Q_p$, is called the norm map. If $x,y\in K$, $\lambda \in \mathbb Q_p$, then $N(xy)=N(x)N(y)$, $N(\lambda x)=\lambda^nN(x)$. The norm map is used to define the normalized absolute value on $K$: $\| x\|=|N(x)|_p$ making $K$ a locally compact totally disconnected topological field.

Denote
$$
O=\{ x\in K:\ \| x\| \le 1\},\ P=\{ x\in K:\ \| x\| <1\},\ U=O\setminus P.
$$
$O$ is a subring of $K$ called the ring of integers, $P$ is an ideal in $O$ called the prime ideal. The multiplicative subgroup $U$ is called the group of units. The quotient $O/P$ is a finite field of characteristic $p$ consisting of $q=p^\upsilon$ elements ($\upsilon \in \mathbb N$). The field $O/P$ is isomorphic to the standard finite field $\mathbb F_q$ consisting of $q$ elements (see \cite{LN}). The normalized absolute value $\| \cdot \|$ takes the values $q^N$, $N\in \mathbb Z$, and 0.

An extension $K$ of degree $n$ is called unramified, if $\| p\| =q^{-1}$. In this case, $P=pO$, $q=p^n$. It is known that, for any $n\in \mathbb N$, there exists an unramified extension of $\mathbb Q_p$ of degree $n$; it is unique up to an isomorphism. Below we fix $n$ and reserve the letter $K$ for this extension. It is generated over $\mathbb Q_p$ by a primitive root of 1 of degree $q-1$. Thus, $K$ contains the group $\mu_{q-1}$ of all the roots of 1 of this degree. On the other hand, if $p\ne 2$, then $K$ does not contain nontrivial roots of 1 of degree p. The Galois group of the extension $K$, that is the group of automorphisms of the field $K$ fixing $Q_p$, is a cyclic group generated by the Frobenius automorphism $\mathfrak g$. On $\mu_{q-1}$, $\mathfrak g$ acts by the rule $\mathfrak g (\omega )=\omega^p$, permuting the roots of 1. The norm map is invariant with respect to $\mathfrak g$: $N(\mathfrak g(x))=N(x)$. On the finite field $O/P$, the Galois group induces the Galois group of $\mathbb F_q$ over $\mathbb F_p$; the automorphism $\mathfrak g$ turns into its finite field counterpart given by raising to the power $p$.

Let $\theta_1,\ldots ,\theta_n\in O$ be such elements that their images in $O/P$ form a basis in $O/P$ over $\mathbb F_p$. Then $\theta_1,\ldots ,\theta_n$ form a basis of $K$ over $\mathbb Q_p$ (called a canonical basis). The choice of this basis determines an identification of $K$ and $\mathbb Q_p^n$. The isomorphism $\mathbb Q_p^n \to K$ as vector spaces over $\mathbb Q_p$ defining this identification has the form
$$
(x_1,\ldots ,x_n)\mapsto \sum \limits_{j=1}^nx_j\theta_j,\quad x_1,\ldots ,x_n\in \mathbb Q_p.
$$
The normalized absolute value on $K$ has the following expression: if $x=\sum \limits_{j=1}^nx_j\theta_j$, then
\begin{equation}
\| x\| =\left( \max\limits_{1\le j\le n}|x_j|_p\right)^n
\end{equation}
(to avoid confusion, note that here and below we consider only the unramified extensions). Below, it will be convenient to assume that $\theta_n=1$.

The multiplicative group $K^*=K\setminus \{ 0\}$ can be described as follows (we consider only the case where $p\ne 2$). If $x\in K^*$, then
\begin{equation}
x=p^\nu \omega \left\{ \prod\limits_{j=1}^{n-1}(1+\theta_jp)^{b_j}\right\} (1+p)^{b_n}
\end{equation}
where $\nu \in \mathbb Z$, $\omega \in \mu_{q-1}$, $b_j\in \mathbb Z_p$ ($j=1,\ldots ,n$). The expression $(1+z)^\beta$, with $z\in K$, $\| z\| <1$, $\beta \in \mathbb Z_p$, is defined as a limit of $(1+z)^{\beta_m}$ where $\beta_m\in \mathbb N$, $\beta_m\to \beta$, as $m\to \infty$, in the topology of $\mathbb Z_p$. An equivalent definition is via the Mahler expansion
\begin{equation}
(1+z)^\beta =1+\sum\limits_{i=1}^\infty z^i\frac{\beta (\beta -1)\cdots (\beta -i+1)}{i!}\quad (\beta \in \mathbb Z_p,\ \|z\|<1)
\end{equation}
convergent in the Banach space of continuous functions on $\mathbb Z_p$ with values from $K$. Obviously, an element (2) has the absolute value $q^{-\nu }$; all the factors in the right-hand side of (2), except the first one, belong to $U$. The elements $\nu ,\omega ,b_1,\ldots ,b_n$ are determined by $x$ in a unique way (thus, for a fixed $x$, each factor $(1+\theta_jp)^{b_j}$ contains a fixed $p$-adic integer $b_j$, so that $(1+\theta_jp)^{b_j}$ is just an element from $U$).

Another canonical representation of an element $x\in K^*$ is
\begin{equation}
x=p^\nu \omega \left( 1+x_1p+x_2p^2+\cdots \right)
\end{equation}
where the first two factors are the same as in (2), $x_1,x_2,\ldots \in \mu_{q-1}\cup \{0\}$, the series converges in $K$, and all the ingredients of (4) are determined in a unique way.

The set $U_1$ of elements (4) with $\nu =0$ and $\omega =1$ is a multiplicative group called the group of principal units. For the unramified extension considered here, the norm map $N$ maps $U_1$ onto the group $U_1(\mathbb Q_p )$ of principal units of the field of $p$-adic numbers. If $\zeta \in U_1(\mathbb Q_p )$, that is $\zeta =1+\zeta_1p+\zeta_2p^2+\cdots$, $\zeta_j\in \mu_{p-1}\cup \{0\}$, the powers $\zeta^\beta$, $\beta \in \mathbb Z_p$, are defined in accordance with (3) and belong to $U_1(\mathbb Q_p )$.

In particular, if $p$ does not divide $n$, then $\frac{1}n\in \mathbb Z_p$, and we have a well-defined root $\zeta^{1/n}\in U_1(\mathbb Q_p )$. Thus, in this case, for any $x\in K^*$ of the form (4), we may write
$$
N(\omega^{-1}x)=p^{n\nu }N(1+x_1p+x_2p^2+\cdots )
$$
and define
\begin{equation}
r=(N(\omega^{-1}x))^{1/n}=p^\nu (N(1+x_1p+x_2p^2+\cdots ))^{1/n}
\end{equation}
as an element of
\begin{equation}
\mathbb Q_p^{(1)} =\left\{ \zeta \in \mathbb Q_p :\ \zeta =p^\nu \left(1+\zeta_1p+\zeta_2p^2+\cdots \right), \ \nu \in \mathbb Z,\ \zeta_j\in \mu_{p-1}\cup \{0\}\right\},
\end{equation}
a multiplicative subgroup of $\mathbb Q_p$.

\section{Spherical coordinates}

For $x\in K^*$, we consider the following elements. Let $\omega =\omega (x)\in \mu_{q-1}$ be the element from (2) or (4). If $p$ does not divide $n$, set
$$
r=r(x)=(N(\omega^{-1}x))^{1/n}\in \mathbb Q_p^{(1)} .
$$
Finally, let
$$
\xi =\xi (x)=\omega^{-1}r^{-1}x,
$$
so that
\begin{equation}
x=\omega (x)\xi (x)r(x).
\end{equation}
We call $(\omega ,\xi ,r)$ {\it the spherical coordinates} of an element $x\in K^*$.

Denote by $\Sigma_n$ the compact multiplicative group
$$
\Sigma_n=\left\{ y\in K^*:\ \omega (y)=1,\ N(y)=1\right\} .
$$

\medskip
\begin{teo}
If $p\ne 2$ and $p$ does not divide $n$, then, for each $x\in K^*$, $\xi (x)\in \Sigma_n$. The representation of an element $x\in K^*$ as a product of elements from $\mu_{q-1}$, $\Sigma_n$, and $\mathbb Q_p^{(1)}$, is unique. The decomposition (7) defines an isomorphism $K^*\cong \mu_{q-1}\times \Sigma_n\times \mathbb Q_p^{(1)}$ of multiplicative topological groups.
\end{teo}

\medskip
{\it Proof}. We will use the representation (2), not with an arbitrary system $\{ \theta_j\}$, but with a special one. In order to construct the latter, we begin with an arbitrary canonical basis $\theta_1,\ldots ,\theta_{n-1},\theta_n$, where $\theta_n=1$. Consider the elements
\begin{equation}
\varepsilon_j=\frac{\theta_j-\mathfrak g(\theta_j)}{1+\mathfrak g(\theta_j)p},\quad j=1,\ldots ,n-1.
\end{equation}
Their images $\overline{\varepsilon_j}$ in $O/P\cong \mathbb F_q$ have the form $\overline{\varepsilon_j}=\overline{\theta_j}-\overline{\theta_j}^p$ where $\overline{\theta_j}$ is the image of $\theta_j$.

Let us show that the elements $\overline{\varepsilon_1},\ldots ,\overline{\varepsilon_{n-1}},1$ form a basis in $\mathbb F_q$ over $\mathbb F_p$. It is sufficient to prove their linear independence.

First we prove the linear independence of $\overline{\varepsilon_1},\ldots ,\overline{\varepsilon_{n-1}}$. Let $c_j\in \mathbb F_p$, $\sum\limits_{j=1}^{n-1}c_j\overline{\varepsilon_j}=0$. Since $c_j^p=c_j$, we have
$$
0=\sum\limits_{j=1}^{n-1}c_j\overline{\theta_j}-
\sum\limits_{j=1}^{n-1}c_j\overline{\theta_j}^p
=\sum\limits_{j=1}^{n-1}c_j\overline{\theta_j}-
\left(\sum\limits_{j=1}^{n-1}c_j\overline{\theta_j}\right)^p,
$$
so that $\sum\limits_{j=1}^{n-1}c_j\overline{\theta_j}
\overset{\text{def}}{=}\lambda \in \mathbb F_p$ and
$$
\sum\limits_{j=1}^{n-1}c_j\overline{\varepsilon_j}-\lambda \overline{\theta_n}=0\quad (\overline{\theta_n}=1).
$$
Since $\overline{\theta_1},\ldots ,\overline{\theta_{n-1}},\overline{\theta_n}$ are linearly independent, we find that $c_1=\ldots =c_{n-1}=\lambda =0$, which proves the linear independence of $\overline{\varepsilon_1},\ldots ,\overline{\varepsilon_{n-1}}$.

Now, let $d_1,\ldots ,d_n\in \mathbb F_p$
\begin{equation}
d_1\left( \overline{\theta_1}-\overline{\theta_1}^p\right)
+d_2\left( \overline{\theta_2}-\overline{\theta_2}^p\right)
+\cdots +d_{n-1}\left( \overline{\theta_{n-1}}-\overline{\theta_{n-1}}^p\right) +d_n=0.
\end{equation}
Raising to the power p we obtain successively that
\begin{gather*}
d_1\left( \overline{\theta_1}^p-\overline{\theta_1}^{p^2}\right)
+d_2\left( \overline{\theta_2}^p-\overline{\theta_2}^{p^2}\right)
+\cdots +d_{n-1}\left( \overline{\theta_{n-1}}^p-\overline{\theta_{n-1}}^{p^2}\right) +d_n=0,\\
...........................................................
...........................................................\\
d_1\left( \overline{\theta_1}^{p^{n-1}}-\overline{\theta_1}^{p^n}\right)
+d_2\left( \overline{\theta_2}^{p^{n-1}}-\overline{\theta_2}^{p^n}\right)
+\cdots +d_{n-1}\left( \overline{\theta_{n-1}}^{p^{n-1}}-\overline{\theta_{n-1}}^{p^n}\right) +d_n=0,\\
d_1\left( \overline{\theta_1}^{p^n}-\overline{\theta_1}^{p^{n+1}}\right)
+d_2\left( \overline{\theta_2}^{p^n}-\overline{\theta_2}^{p^{n+1}}\right)
+\cdots +d_{n-1}\left( \overline{\theta_{n-1}}^{p^n}-\overline{\theta_{n-1}}^{p^{n+1}}\right) +d_n=0.
\end{gather*}
Note that $\overline{\theta_j}^{p^{n+1}}=\overline{\theta_j}^{qp}=\overline{\theta_j}^p$
and add up all the equalities. We find that $nd_n=0$ in $\mathbb F_p$, and since $p$ does not divide $n$, $d_n=0$. Then (9) implies the equalities $d_1=\ldots =d_{n-1}=0$.

Thus, we know that the collection $\overline{\varepsilon_1},\ldots ,\overline{\varepsilon_{n-1}},1$ is a $\mathbb F_p$-basis in $O/P$. Therefore we may write a representation like (2) based on the canonical basis $\varepsilon_1,\ldots ,\varepsilon_{n-1},1$: for every $x\in K^*$,
\begin{equation}
x=p^\nu \omega \left\{ \prod\limits_{j=1}^{n-1}(1+\varepsilon_jp)^{b_j'}\right\} (1+p)^{b_n'}
\end{equation}
where $\nu \in \mathbb Z$, $\omega =\omega (x)\in \mu_{q-1}$, $b_1',\ldots ,b_n'\in \mathbb Z_p$.

By (8), we can write
$$
1+\varepsilon_jp=\frac{1+\theta_jp}{1+\mathfrak g(\theta_j)p},
$$
whence
\begin{equation}
\prod\limits_{j=1}^{n-1}(1+\varepsilon_jp)^{b_j'}=\frac{y}{\mathfrak g(y)},\quad y=\prod\limits_{j=1}^{n-1}(1+\theta_jp)^{b_j'}.
\end{equation}
It follows from (11) that
$$
\prod\limits_{j=1}^{n-1}(1+\varepsilon_jp)^{b_j'}\in \Sigma_n.
$$

On the other hand, $p^\nu (1+p)^{b_n'}\in \mathbb Q_p^{(1)}$, and $N(\omega^{-1}x)=p^{n\nu }(1+p)^{nb_n'}$. Thus
$$
r(x)=p^\nu (1+p)^{b_n'},
$$
and we have got the representation (7) with $\xi (x)=\dfrac{y}{\mathfrak g(y)}$.

If we have another representation $x=\omega_1\xi_1r_1$, $\omega_1\in \mu_{q-1}$, $\xi_1\in \Sigma_n$, and $r_1\in \mathbb Q_p^{(1)}$, then it follows directly from the definitions of  $\Sigma_n$ and $\mathbb Q_p^{(1)}$ that $\omega_1=\omega$. Then, applying $N$ we get $N(\omega^{-1}x)=r_1^n$, so that $r_1=r$ and $\xi_1=\xi$.

The fact that the representation (7) is compatible with the algebraic operations and the topologies on the corresponding topological groups follows immediately from the properties of the norm map and the group of principal units as a topological $\mathbb Z_p$-module. $\qquad \blacksquare$

\medskip
Below we will often denote our spherical coordinates by $x=(\eta ,r)$ where $\eta =(\omega ,\xi )\in Z_n=\mu_{q-1}\times \Sigma_n$. As in the classical situation, sometimes it is convenient to extend the spherical coordinates to the whole of $K$ -- for $x=0$, we set $r=0$ while $\eta$ is not defined.

In order to derive a formula for integration in spherical coordinates, denote by $d\xi$ the Haar measure on $\Sigma_n$ normalized by the relation $\int_{\Sigma_n}d\xi =1$. A Haar measure on $\mathbb Q_p^{(1)}$ is induced by the multiplicative Haar measure on $\mathbb Q_p^*$ having the form $\dfrac{dr}{|r|_p}$ where $dr$ is the additive Haar measure. Similarly (see \cite{B}), $\dfrac{dx}{\|x\|}$ is the Haar measure on $K^*$. As usual, we assume that $\int_O dx=\int_{\mathbb Z_p}dr=1$.

For the direct product $K^*=\mu_{q-1}\times \Sigma_n\times \mathbb Q_p^{(1)}$, we have the integration formula
\begin{equation}
\int\limits_{K^*}f(x)\frac{dx}{\|x\|}=c\sum\limits_{\omega \in \mu_{q-1}}\int\limits_{\Sigma_n}d\xi \int\limits_{\mathbb Q_p^{(1)}}f(\omega \xi r)\frac{dr}{|r|_p}
\end{equation}
valid, for example, for any continuous function on $K^*$ with a compact support. In order to find the normalization constant $c$, we take for $f$ the indicator function of the group of units $U$.

It is known \cite{K} that
$$
\int\limits_Udx=1-q^{-1},
$$
$$
\int\limits_{r\in \mathbb Q_p^{(1)} ,|r|_p=1}dr=\frac{1}{p-1}(1-\frac{1}p)=p^{-1}.
$$
Therefore $1-q^{-1}=c(q-1)p^{-1}$, whence $c=\dfrac{1}{p^{n-1}}$.

It is easy to rewrite (12) in terms of additive Haar measures. Substituting $f(x)\|x\|$ for $f(x)$ in (12) we find that
\begin{equation}
\int\limits_K f(x)\,dx=\frac{1}{p^{n-1}}\sum\limits_{\omega \in \mu_{q-1}}\int\limits_{\Sigma_n}d\xi \int\limits_{\mathbb Q_p^{(1)}}f(\omega \xi r)|r|_p^{n-1}\,dr.
\end{equation}
We will not study exact conditions on $f$, under which (13) is valid. It is sufficient if $f$ is continuous on $K$ and has a compact support.

\section{Homogeneous Distributions}

Let $\pi :\ \mathbb Q_p^* \to \mathbb C$ be a multiplicative quasicharacter, that is $\pi (z)=|z|_p^s\theta (z)$ where $s\in \mathbb R$, $\theta (z)$ is a multiplicative character, $|\theta (z)|=1$, such that $\theta (p)=1$.

A continuous function $f:\ \mathbb Q_p^n \to \mathbb C$ is called a homogeneous function of degree $\pi$, if
\begin{equation}
f(\lambda x_1,\ldots ,\lambda x_n)=\pi (\lambda )f(x_1,\ldots ,x_n)
\end{equation}
for any $\lambda \in \mathbb Q_p^{(1)}$, $x=(x_1,\ldots ,x_n)\in \mathbb Q_p^n$.

A Bruhat-Schwartz distribution $f\in \mathcal D'(\mathbb Q_p^n )$ is called a homogeneous distribution of degree $\pi$, if
\begin{equation}
\langle f,\varphi_\lambda \rangle =\pi (\lambda )|\lambda |_p^n\langle f,\varphi \rangle ,\quad \varphi_\lambda (x)=\varphi (\lambda^{-1}x)\ (x\in \mathbb Q_p^n ),
\end{equation}
for any $\varphi \in \mathcal D(\mathbb Q_p^n )$, $\lambda \in \mathbb Q_p^{(1)}$. The definitions (14) and (15) are slightly more general than the usual ones \cite{GGP,VVZ,AKS} -- we take only ``positive'' $\lambda \in \mathbb Q_p^{(1)}$.

As before, studying the structure of homogeneous distributions we identify $\mathbb Q_p^n$ with the unramified extension $K$. The definition (15) makes sense in this case too. {\bf Below we assume that} $p\ne 2$ {\bf and} $p$ {\bf does not divide} $n$. Then we may use the spherical coordinates $x=\omega \xi r$ ($x\in K^*$, $\omega \in \mu_{q-1}$, $\xi \in \Sigma_n$, and $r\in \mathbb Q_p^{(1)}$).

It follows from the representations (10) and (11) that the group $\Sigma_n$ is isomorphic to the direct product of $n-1$ copies of $\mathbb Z_p$. Therefore we have natural spaces of test functions $\mathcal D(\Sigma_n)$ (consisting of locally constant functions on $\Sigma_n\cong \mathbb Z_p^{n-1}$) and $\mathcal D(Z_n)$, $Z_n=\mu_{q-1}\times \Sigma_n$, as well as the spaces of distributions $\mathcal D'(\Sigma_n)$ and $\mathcal D'(Z_n)$. If $F\in \mathcal D'(Z_n)$ is generated by an ordinary function $F(\omega ,\xi )$, that means that
$$
\langle F,\psi \rangle =\frac{1}{q-1}\sum\limits_{\omega \in \mu_{q-1}}\int\limits_{\Sigma_n}F(\omega ,\xi )\psi (\omega ,\xi )\,d\xi ,\quad \psi \in \mathcal D(Z_n).
$$

If $f$ is a continuous homogeneous function, then it follows from (14) that
\begin{equation}
f(\omega \xi r)=\pi (r)f(\omega \xi ).
\end{equation}

A function $f$ of the form (16) with $\pi (r)=|r|_p^s\theta (r)$, $\mathrm{Re}\ s>-n$, determines a distribution from $\mathcal D'(K)$ in a straightforward way. Using the integration formula (13) we get, for any $\varphi \in \mathcal D(K)$, that
$$
\langle f,\varphi \rangle =\int\limits_K f(x)\varphi (x)\,dx\\
=p^{1-n}\sum\limits_{\omega \in \mu_{q-1}} \int\limits_{\mathbb Q_p^{(1)}}|r|_p^{s+n-1}\theta (r)\,dr\int\limits_{\Sigma_n}f(\omega \xi )\varphi (r\omega \xi )\,d\xi.
$$

More generally, if $\mathrm{Re}\ s>-n$, $F\in \mathcal D'(Z_n)$, then the distribution $f=\pi (r)F$ is given by the relation
$$
\langle f,\varphi \rangle =p^{1-n}(p^n-1) \int\limits_{\mathbb Q_p^{(1)}}\langle F,\varphi (r\cdot )\rangle |r|_p^{s+n-1}\theta (r)\,dr,
$$
$\varphi \in \mathcal D(K)$.

The function $r\mapsto \langle F,\varphi (r\cdot )\rangle$ is locally constant and has a compact support by virtue of the compactness of $Z_n$. In particular, suppose that $\varphi (r\omega \xi )=0$ if $|r|_p>p^\nu$. Then
\begin{multline}
\langle f,\varphi \rangle =p^{1-n}\sum\limits_{\omega \in \mu_{q-1}} \left\{ \int\limits_{r\in \mathbb Q_p^{(1)} :\ |r|_p\le p^\nu }\langle F,\varphi (r\cdot )-\varphi (0)\rangle |r|_p^{s+n-1}\theta (r)\,dr\right. \\ \left.
-\varphi (0)\langle F,1\rangle \int\limits_{r\in \mathbb Q_p^{(1)} :\ |r|_p\le p^\nu }|r|_p^{s+n-1}\theta (r)\,dr\right\}
\end{multline}

The first integral in (17) is an entire function of $s$. Next,
\begin{equation}
\int\limits_{r\in \mathbb Q_p^{(1)} :\ |r|_p\le p^\nu }|r|_p^{s+n-1}\theta (r)\,dr=\sum\limits_{j=-\infty }^\nu p^{j(s+n-1)}\int\limits_{r\in \mathbb Q_p^{(1)} :\ |r|_p=p^j}\theta (r)\,dr
\end{equation}

If the character $\theta$ is nontrivial on the group of principal units of $\mathbb Q_p$, then (making the change of variables $r=a\rho$ where $a$ is a principal unit with $\theta (a)\ne 1$) we find that all the integrals in the right-hand side of (18) equal zero. If $\theta (r)\equiv 1$, then \cite{V}
$$
\int\limits_{r\in \mathbb Q_p^{(1)} :\ |r|_p=p^j}dr=p^{j-1},
$$
so that
$$
\int\limits_{\mathbb Q_p^{(1)} }|r|_p^{s+n-1}\,dr=p^{-1}\sum\limits_{j=-\infty }^\nu p^{j(s+n)}=\frac{p^{\nu (s+n)-1}}{1-p^{-s-n}}.
$$

It follows that the distribution $f=\pi (r)F$ from $\mathcal D'(K)$ is defined by the analytic continuation procedure for any quasicharacter $\pi$ and any distribution $F\in \mathcal D'(Z_n)$, with a single exception, the case where $\pi (r)=|r|_p^{-n}$ (this quasicharacter will be called exceptional) and $\langle F,1\rangle \ne 0$. The residue of $\langle f,\varphi \rangle$ at $s=-n+\frac{2\pi ki}{\log p}$ equals
$$
p^{1-n}\varphi (0)\langle F,1\rangle \mathrm{Res}_{s=-n}
\frac{p^{\nu (s+n)-1}}{1-p^{-s-n}}=\frac{1}{p^n\log p}\varphi (0)\langle F,1\rangle ,
$$
so that
$$
\mathrm{Res}_{s=0}f=\frac{\langle F,1\rangle }{p^n\log p}\delta .
$$

For a non-exceptional quasicharacter $\pi$, the distribution $f=\pi (r)F$ obviously satisfies (15). Below we show that the above construction covers the whole class of homogeneous distributions. However we need some auxiliary results.

Denote $\mathbb Q_p^{(+)} =\mathbb Q_p^{(1)} \cup \{ 0\}$. With the metric, topology, and (additive) Haar measure induced on $\mathbb Q_p^{(+)}$ from $\mathbb Q_p$, we denote by $\mathcal D(\mathbb Q_p^{(+)} )$ the space of all locally constant functions on $\mathbb Q_p^{(+)}$ with compact supports, and by $\mathcal D'(\mathbb Q_p^{(+)} )$ the dual space containing in a standard manner all the ``ordinary'' functions. The definition (15) (with $n=1$) of a homogeneous distribution makes sense for distributions on $\mathbb Q_p^{(+)}$.

\medskip
\begin{lem}
A homogeneous distribution on $\mathbb Q_p^{(+)}$ of degree $\pi$, where $\pi (r)\ne |r|_p^{-1}$, has the form $C\pi$, $C=\text{const}$.
\end{lem}

\medskip
The {\it proof} is identical to the one known for distributions on $\mathbb Q_p$ (see \cite{VVZ}).

Consider an arbitrary test function (a locally constant function with a compact support) $\varphi \in \mathcal D(K)$.

\medskip
\begin{lem}
The function $\varphi$ admits a decomposition into a finite sum
\begin{equation}
\varphi (\omega \xi r)=\varphi (0)\Delta_l(r)+\sum\limits_{m=1}^M\varphi (\omega \xi r_m)\Delta_l(r-r_m),\quad \omega \in \mu_{q-1}, \xi \in \Sigma_n, r\in \mathbb Q_p^{(+)} ,
\end{equation}
where $r_m$ are some points of $\mathbb Q_p^{(1)}$ depending only on the function $\varphi$, $\Delta_l(z)$ is the indicator function of some ball $\{ z\in \mathbb Q_p :\ |z|_p\le p^l\}$ ($l\in \mathbb Z$).
\end{lem}

\medskip
{\it Proof}. Suppose first that $\varphi (0)=0$, that is $\text{supp } \varphi \subset C$ where
$$
C=\left\{ x\in K:\ q^\nu \le \|x\| \le q^N\right\},\quad \nu \le N\ (\nu ,N\in \mathbb Z),
$$
and $\varphi (x+y)=\varphi (x)$ for any $x\in K$, if $\|y\| \le q^l$, and we may assume that $l<\nu$.

In spherical coordinates, we have $C=\mu_{q-1}\times \Sigma_n\times \tilde{C}$ where
$$
\tilde{C}=\left\{ r\in \mathbb Q_p^{(1)} :\ p^\nu \le |r|_p \le p^N\right\} ,
$$
and the above local constancy of $\varphi$ is equivalent to the local constancy of the function $r\mapsto \varphi (\omega \xi r)$:
$$
\varphi (\omega \xi (r+r'))=\varphi (\omega \xi r),\quad \text{if $|r'|_p\le p^l$.}
$$

Let us take a finite covering of $\tilde{C}$ by non-intersecting balls of radius $p^l$ with the centers $r_m$, $m=1,\ldots ,M$. We get the representation
\begin{equation}
\varphi (\omega \xi r)=\sum\limits_{m=1}^M\varphi (\omega \xi r_m)\Delta_l(r-r_m).
\end{equation}

If $\varphi (0)\ne 0$, consider the function $\varphi_1 (\omega \xi r)=\varphi (\omega \xi r)-\varphi (0)\Delta_l(r)$. Clearly, $\varphi_1\in \mathcal D(K)$,
$$
\varphi_1(x)=\begin{cases}
0, & \text{if $\|x\|\le q^l$};\\
\varphi (x), & \text{if $\|x\|>q^l$}.
\end{cases}
$$
Applying (20) to the function $\varphi_1$ and noticing that $\Delta_l(r)\sum\limits_{m=1}^M\Delta_l(r-r_m)\equiv 0$, we obtain the equality (19). $\qquad \blacksquare$

\medskip
Now we can give a description of homogeneous distributions in the sense of (15).

\medskip
\begin{teo}
Suppose that $p\ne 2$, $p$ does not divide $n$, and $\pi$ is a non-exceptional quasicharacter. Then any homogeneous distribution $f$ of degree $\pi$ has the form $f=\pi (r)F$, $F\in \mathcal D'(Z_n)$.
\end{teo}

\medskip
{\it Proof}. Let $\varphi \in \mathcal D(\mathbb Q_p^{(+)} )$, $\varphi (0)=0$, $\psi \in \mathcal D(Z_n)$. Then $(\varphi \otimes \psi )(\omega \xi r)=\psi (\omega \xi )\varphi (r)$ is a test function from $\mathcal D(K)$. The linear mapping
$$
\varphi \mapsto \langle f,\varphi \otimes \psi \rangle
$$
is a homogeneous distribution from $\mathcal D'(\mathbb Q_p^{(+)} )$ of degree $\pi_1$, $\pi_1(\lambda )=\pi (\lambda )|\lambda |_p^{n-1}$. By Lemma 1, for each $\varphi$,
$$
\langle f,\varphi \otimes \psi \rangle =C_\psi \langle \pi_1,\varphi \rangle
$$
where $C_\psi$ is some constant.

Let $\varphi \in \mathcal D(\mathbb Q_p^{(+)} )$ be such a test function that $\varphi (0)=0$ and $\langle \pi_1,\varphi \rangle =\dfrac{p^{n-1}}{p^n-1}$. Define a distribution $F\in \mathcal D'(Z_n)$ setting
$$
\langle F,\psi \rangle =\langle f,\varphi \otimes \psi \rangle ,\quad \psi \in \mathcal D(Z_n).
$$
For any $\Phi \in \mathcal D(\mathbb Q_p^{(+)} )$, such that $\Phi (0)=0$, we have
\begin{multline*}
\langle \pi F-f,\Phi \otimes \psi \rangle =p^{1-n}(p^n-1)\int\limits_{\mathbb Q_p^{(+)} }\pi_1(r)\langle F,\psi \rangle\Phi (r)\,dr-\langle f,\Phi \otimes \psi \rangle \\
=p^{1-n}(p^n-1)\langle F,\psi \rangle \langle \pi_1 ,\Phi \rangle
=C_\psi p^{1-n}(p^n-1)\langle \pi_1,\varphi \rangle \langle \pi_1,\Phi \rangle -C_\psi \langle \pi_1,\Phi \rangle =0.
\end{multline*}

Taking into account Lemma 2, we see that the distribution $\pi F-f\in \mathcal D'(K)$ is concentrated at the origin. Therefore \cite{GGP,VVZ} there exists such a constant $c\in \mathbb C$ that  $\pi F-f=c\delta$. Since  $\pi F-f$ is a homogeneous distribution of degree $\pi$, and $\pi$ is non-exceptional, we find that $c=0$, as desired. $\qquad \blacksquare$

\section{Skew Product Decompositions of $p$-Adic L\'evy Processes}

As before, we assume that $p\ne 2$ and $p$ does not divide $n$, and identify $\mathbb Q_p^n$ with the unramified extension $K$.

Let $X_t$ be a rotation-invariant termporally and spatially homogeneous process on $K$ with independent increments. Its transition probability has the form
$$
P_t(x,B)=\Phi_t(B-x),\quad x\in K,\ B\in \mathcal B(K),
$$
where $\Phi_t$ is a semigroup of measures, $\mathcal B(K)$
denotes the Borel $\sigma$-algebra of $K$. Below we assume that $\Phi_t$ is absolutely continuous with respect to the Haar measure (see \cite{K,Y} for a complete description of such processes). In this case
$$
P_t(x,B)=\int\limits_B\Gamma (t, \|x-y\|)\,dy,\quad x\in K,\ B\in \mathcal B(K).
$$
See \cite{K,Y} for an explicit expression of the density $\Gamma$.

Consider the process $R_t=r(X_t)$.

\medskip
\begin{teo}
The process $R_t$ is a Markov process in $\mathbb Q_p^{(+)}$ with the transition probability
$$
\tilde{P}_t(r(x),\tilde{B})=P_t(x,r^{-1}(\tilde{B})),\quad x\in K,\ \tilde{B}\in \mathcal B(\mathbb Q_p^{(+)} ).
$$
\end{teo}

\medskip
{\it Proof}. By the general theorem on transformations of the phase space of a Markov process (\cite{Dyn}, Theorem 10.13), it is sufficient to check that if $x'\in K$, $r(x')=r(x)$, then
\begin{equation}
P_t(x,r^{-1}(\tilde{B}))=P_t(x',r^{-1}(\tilde{B})).
\end{equation}

For the above $x$ and $x'$, we have $\|x\|=|N(x)|_p$, so that $\|x\|=\|\omega^{-1}(x)x\|=|N(\omega^{-1}x)|_p=|r(x)|^n_p$, and $\|x'\|=\|x\|$. On the other hand,
$$
P_t(x,r^{-1}(\tilde{B}))=\int\limits_{r(y)\in \tilde{B}}\Gamma (t,\|x\|\cdot \|1-x^{-1}y\|)\,dy.
$$
The change of variables $x^{-1}y=z$ yields
\begin{multline*}
P_t(x,r^{-1}(\tilde{B}))=\|x\| \int\limits_{r(x)r(z)\in \tilde{B}}\Gamma (t,\|x\|\cdot \|1-z\|)\,dz\\
=\|x'\| \int\limits_{r(x')r(z)\in \tilde{B}}\Gamma (t,\|x'\|\cdot \|1-z\|)\,dz=P_t(x',r^{-1}(\tilde{B})),
\end{multline*}
and we have proved the equality (21).$\qquad \blacksquare$

\medskip
Let us turn to the ``angular'' process $z_t=\eta (X_t)$ on $Z_n=\mu_{q-1}\times \Sigma_n$ generated by the process $X_t$. Let $D_T(K)$ be the space of c\`adl\`ag functions $[0,T]\to K$ endowed with the Skorokhod topology (see \cite{GS} regarding this notion for functions with values from a metric space). The mappings $r$ and $\eta$ induce the corresponding mappings $D_T(K)\to D_T(\mathbb Q_p^{(1)} )$ and $D_T(K)\to D_T(Z_n)$. Let $\mathbf P_x$ be the probability measure on $D_T(K)$ corresponding to the process $X_t$ with the starting point $x$. Denote by $\zeta$ the first hitting time for the point 0, that is
$$
\zeta =\inf \{ t>0:\ X_t=0\text{ or $X_{t-}=0$}\}
$$
where we assume that $\inf \varnothing =\infty$.

Let $\mathcal F_{0,T}^{\mathbb Q_p^{(+)}}$ be the $\sigma$-algebra on $D_T(\mathbb Q_p^{(+)} )$ generated by the process $R_t$. It induces the $\sigma$-algebra $r^{-1}\left( \mathcal F_{0,T}^{\mathbb Q_p^{(+)}} \right)$ on $D_T(K)$. By the existence of regular conditional distributions (see Theorem 5.3 in \cite{Kall}), there exists a probability kernel $W_z^{R(\cdot )}$ from $Z_n\times D_T(\mathbb Q_p^{(1)} )$ to $D_T(Z_n)$, such that for any $x\in K$, $x\ne 0$, and any measurable $F\subset D_T(Z_n)$,
$$
W_{\eta (x)}^{r(X(\cdot ))}=\mathbf P_x\Biggl[ X(\cdot )\in \eta^{-1}(F)\Biggr| r^{-1}\left(\mathcal F_{0,T}^{\mathbb Q_p^{(+)}} \right) \Biggr] ,
$$
for $\mathbf P_x$-almost all $X(\cdot )$ in $[\zeta >T]\subset D_T(K)$.

The probability measure $W_z^{R(\cdot )}$ is considered as the conditional distribution of the process $z_t$, given $z_0=z$ and a radial path $R(\cdot )$ in $D_T(\mathbb Q_p^{(1)} )$.

It can be proved in exactly the same way as in \cite{L} that for $x\ne 0$, $r\mathbf P_x$-almost all $R(\cdot )$ in $[\zeta >T]\subset D_T(\mathbb Q_p^{(1)} )$, the process $z_t$ is a non-homogeneous L\'evy process, that is there exists a two-parameter semigroup of random measures $\nu_{s,t}$ on $Z_n$, such that for any natural number $m$, any continuous function $f$ on $Z_n^m$, any points $t_1<t_2<\ldots <t_m$ from $[0,T]$, and for $x=\rho z$, $\rho \in \mathbb Q_p^{(1)}$, $z\in Z_n$,
\begin{multline*}
\mathbf E_{\rho z}\Biggl[ f(z_{t_1},\ldots ,z_{t_m})\Biggr| \mathcal F_{0,T}^{\mathbb Q_p^{(+)}} \Biggr] \\
=\int\limits_{Z_n^m}f(zz_1,zz_1z_2,\ldots ,zz_1\cdots z_m)\nu_{0,t_1}(dz_1)\nu_{t_2,t_1}(dz_2)\cdots \nu_{t_m,t_{m-1}}(dz_m).
\end{multline*}

As in the classical situation, the pair $(R_t,z_t)$ forms a ``skew product'' representation of the process $X_t$.


\medskip
\begin{thebibliography}{999}
\bibitem{AKS}
S. Albeverio, A. Yu. Khrennikov, and V. M. Shelkovich, Harmonic analysis in the $p$-adic Lizorkin spaces: fractional operators, pseudo-differential equations, $p$-adic wavelets, Tauberian
theorems, {\it J. Fourier Anal. Appl.} {\bf 12} (2006), 393--425.
\bibitem{B}
N. Bourbaki, {\it Integration II}, Springer, Berlin, 2004.
\bibitem{Dyn}
E. B. Dynkin, {\it Markov Processes I}, Springer, Berlin, 1965.
\bibitem{FV}
I. B. Fesenko and S. V. Vostokov, {\it Local Fields and Their Extensions},
American Mathematical Society, Providence, 2002.
\bibitem{G}
A. R. Galmarino, Representation of an isotropic diffusion as a skew product, {\it Z. Wahrscheinlichkeitstheor. Verw. Geb.} {\bf 1} (1963), 359--378.
\bibitem{GGP}
I. M. Gelfand, M. I. Graev, and I. I. Piatetski-Shapiro, {\it Representation Theory and Automorphic Functions}, Saunders, Philadelphia, 1969.
\bibitem{GS}
I. I. Gihman and A. V. Skorohod, {\it The Theory of Stochastic
Processes I}, Springer, Berlin, 1974.
\bibitem{Kall}
O. Kallenberg, {\it Foundations of Modern Probability}, Springer, Berlin, 1997.
\bibitem{Koch}
H. Koch, {\it Algebraic Number Fields}, Encycl. Math. Sci. Vol. 62, Springer, Berlin, 1992.
\bibitem{K}
A. N. Kochubei, {\it Pseudo-Differential Equations and Stochastics over Non-Archimedean Fields}, Marcel Dekker, New York, 2001.
\bibitem{Le}
C. Lemoine, Fourier transforms of homogeneous distributions, {\it Ann. Scuola Norm. Sup. Pisa}, {\bf 26} (1972), 117--149.
\bibitem{L}
M. Liao, A decomposition of Markov processes via group actions, {\it J. Theor. Probab.} {\bf 22} (2009), 164--185.
\bibitem{LN}
R. Lidl and H. Niederreiter, {\it Finite Fields}, Cambridge
University Press, 1983.
\bibitem{PR}
E. J. Pauwels and L. C. G. Rogers, Skew-product decompositions of Brownian motions, {\it Contemp. Math.} {\bf 73} (1988), 237--262.
\bibitem{R}
A. M. Robert, {\it A Course in $p$-Adic Analysis}, Springer, New York, 2000.
\bibitem{Sch}
W. Schikhof, {\it Ultrametric Calculus}, Cambridge University
Press, 1984.
\bibitem{V}
V. S. Vladimirov, Generalized functions over the field of $p$-adic numbers,
{\it Russian Math. Surveys}, {\bf 43}, No. 5 (1988), 19--64.
\bibitem{VVZ}
V. S. Vladimirov, I. V. Volovich and E. I. Zelenov, {\it $p$-Adic Analysis and
Mathematical Physics}, World Scientific, Singapore, 1994.
\bibitem{W}
A. Weil, {\it Basic Number Theory}, Springer, Berlin, 1967.
\bibitem{Y}
K. Yasuda, Additive processes on local fields, {\it J. Math. Sci. Univ. Tokyo}, {\bf 3} (1996), 629--654.


\end{thebibliography}
\end{document}